\documentclass[12pt,a4paper,twoside]{amsart}
\usepackage[english]{babel}
\usepackage{amsmath}
\usepackage{amssymb}
\usepackage[latin1]{inputenc}  %<--- fuer UNIXe
\theoremstyle{plain}
\newtheorem{theorem}{Proposition}[section]
\newtheorem{thm}[theorem]{Theorem}
\newtheorem{lemma}[theorem]{Lemma}

\renewcommand{\H}{\mathcal H}
\newcommand{\N}{\mathbb N}

\newcommand{\C}{\mathbb C}

\newcommand{\D}{\mathcal D}
\newcommand{\F}{\mathcal F}

\newcommand{\la}{\big\langle}
\newcommand{\ra}{\big\rangle}
\renewcommand{\L}{\mathcal L^+(\D)}
\begin{document}
\title{Additive Mappings in Algebras of Unbounded Operators Preserving Operators of Rank One}
\author{WERNER TIMMERMANN}

\address{Institut f\"ur Analysis\\
         Technische Universit\"at Dresden\\
         D-01062 Dresden, Germany}
\email{timmerma@math.tu-dresden.de}

\begin{abstract}
Let $\D$ be a dense linear manifold in a Hilbert space $\H$ and let $\L$ be the *-algebra of all linear operators $A$ such that $A\D \subset \D, A^*\D \subset \D$. Denote by $\F(\D)$ the *-ideal of $\L$ consisting of all finite-rank operators. A characterization of the structure of additive mappings on $\F(\D)$ preserving operators of rank one or projections of rank one is given. The corresponding results for algebras of bounded operators on Banach spaces were given in \cite{omladic93}.
\end{abstract}

\thanks{\hspace{-0.5cm} Keywords: algebras of unbounded operators, rank-one preservers  }

\maketitle

\section{Introduction and statement of the results}

In the last 15 years there can be observed a growing interest in so-called linear preserver problems for operator algebras over infinite-dimensional spaces. It appeared that many problems can be reduced to the problem of determining the structure of linear mappings that carry operators of rank one into themselves. The vast literature on linear mappings preserving rank deals mostly with algebras of bounded operators on Banach or Hilbert spaces (see e.g. \cite{ hou89, molnar98d, omladic}  and the references therein). There is only one paper containing results on the structure of mappings preserving rank-one operators in the context of algebras of unbounded operators \cite{ti93}.\\
In many cases it turns out that such mappings are automorphisms or antiautomorphisms of the underlying algebras. On the other hand, it is easy to see that automorphisms and antiautomorphisms preserve operators of rank one as well as idempotents and nilpotent operators.\\
 A more general approach is to consider such algebras only as rings. It was just this point of view that opened the door for a fruitful application of results of abstract ring theory to the investigation of special problems in operator algebras. In particular, there were studied not only linear preservers but also additive ones.\\
In \cite{omladic93} Omladi\u{c} and \u{S}emrl studied additive mappings preserving rank-one operators. Let us fix some notions and notation to formulate some of their results. For a Banach space $\mathcal X$ let us denote by $\mathcal B(\mathcal X)$ the Banach algebra of all bounded operators on $\mathcal X$ and by $\mathcal F(\mathcal X) \subset \mathcal B(\mathcal X)$ the ideal of finite rank operators. $\mathcal X'$ denotes the dual Banach space to $\mathcal X$. For our purpose it is enough to consider infinite dimensional complex Banach spaces. An operator $F \in \mathcal B(\mathcal X)$ is of rank one if and only if $F = f \otimes x, f \in \mathcal X', x \in \mathcal X$, where $Fy = f(y)x, ~ y \in \mathcal X$. Let $\mathcal A$ be a standard operator algebra, i.e. a subalgebra of $ \mathcal B(\mathcal X)$ such that  $\mathcal F(\mathcal X) \subset \mathcal A$. An additive mapping $\Phi\colon \mathcal A \rightarrow \mathcal A $ is said to preserve rank-one operators if $\Phi(A)$ is of rank one if $A$ is so. $\Phi$ preserves rank one operators in both directions if $\Phi(A)$ is of rank one if and only if $A$ is of rank one. $\Phi$ is said to preserve linear spans of projections of rank one if for every projection $P$ of rank one the inclusion $\Phi(\C P) \subset \C \Phi(P)$ holds.\\
Let $h$ be a ring automorphism of $\C$. A mapping $A \colon E \rightarrow F,~ E, F $ - complex vector spaces, is called $h$-quasilinear (conjugate-linear, resp.) if it is additive and satisfies $A(\lambda x) = h(\lambda) Ax$ ($A(\lambda x) = \bar{\lambda} A x$, resp.) for all $\lambda \in \C, x \in E$. \\

{\bf Theorem A}(\cite{omladic93}, Theorem 3.3) \emph{Let X be a complex Banach space, dim $X > 1$. Suppose that $\Phi \colon \mathcal F(\mathcal X) \rightarrow \mathcal F(\mathcal X)$ is an additive mapping preserving operators of rank one in both directions. Then there is a ring automorphism $h: \C \rightarrow \C$ and there are h-quasilinear bijective mappings $A \colon \mathcal X \rightarrow \mathcal X$ and $C \colon \mathcal X' \rightarrow \mathcal X'$ such that
$$\Phi(x \otimes f) = Ax \otimes Cf, \quad x \in \mathcal X, f \in \mathcal X'$$or h-quasilinear bijective mappings $A \colon \mathcal X \rightarrow \mathcal X'$ and $C \colon \mathcal X' \rightarrow \mathcal X$ such that
$$\Phi(x \otimes f) = Cf \otimes Ax, \quad x \in \mathcal X, f \in \mathcal X'$$}

 Theorem 4.4 in \cite{omladic93} is their formulation of the Main Theorem for the case of a Hilbert space. We give this formulation in a slight modification, namely we avoid the explicit use of the Banach space conjugate operator $T'$. Only the Hilbert space adjoint operator $T^*$ is used.\\

{\bf Theorem B}(\cite{omladic93}, Theorem 4.4) \emph{Let $\H$ be an infinite-dimensional complex Hilbert space. Suppose that $\Phi \colon \F(\H) \rightarrow \F(\H)$ is an additive surjective mapping preserving projections of rank one and their linear spans. Then one of the following four possibilities is valid:\\
i) one has
$$\Phi(T) = ATA^{-1}$$
for all $T \in \F(\H)$, where $A \colon \H \rightarrow \H$ is a continuous bijective either linear or conjugate-linear operator;\\
or\\
ii) one has 
$$\Phi(T) = CT^*C^{-1}$$
for all $T \in \F(\H)$, where $C \colon \H \rightarrow \H$ is a continuous bijective either linear or conjugate-linear operator}\\

The proof of Theorem B splits into two parts. The first part is purely algebraic and does not use that $\mathcal X$ (in our formulation $\H$) is a Banach space. The second, analytical part consists in proving the continuity of the representing operators and the continuity of $h$. Here it is used that $\mathcal X$ is a Banach space.\\
Let us remark that Moln\'ar used Theorem B in \cite{molnar00a} to prove a Wigner-type theorem in Banach spaces which also gives a short and elegant proof of Wigner's celebrated unitary-antiunitary theorem \cite{wigner31}.\\ 
The aim of the present paper is to indicate how the above mentioned results fit in the concept of algebras of unbounded operators in Hilbert spaces. In the first part we give a rather general formulation of Theorem A using dual pairs. Then Theorem B is given for algebras of unbounded operators. May be this variant is not the most general and elegant one but it is near to the interests of the author. In a forthcoming paper we will discuss an application of these results to a proof of Wigner-type theorems in the context of algebras of unbounded operators, quite similar to \cite{molnar00a}.\\
Let us fix the necessary notation. A standard reference for algebras of unbounded operators is \cite{schm90}.\\ 
Let $\mathcal D$ be a dense linear manifold in a Hilbert space $\mathcal H$ with scalar
product $\la ~,~ \ra $ (which is supposed to be conjugate linear in the first and linear in the second component). The set of linear operators $\L = \{ A: A \D \subset \D, A^* \D \subset \D\}$ is a $\ast$-algebra with respect to the natural operations and the involution $A \rightarrow A^+ = A^*|\D$.
The graph topology $t$ on $\D$ induced by $\L$ is generated by the directed family of seminorms $\phi \rightarrow ||\phi ||_A = || A \phi ||,~ \forall ~  A \in \L, ~  \phi \in \D$. $\D$ is called an (F)-domain, if $(\D, t)$ is an (F)-space. Remark that in this case the graph topology $t$ can be given by a system of seminorms $\{ \| \cdot \|_n = \|A_n \cdot \|, n \in \N, A_n \in \L\}$ with: 
$$A_1 = I, \  A_n = A^+_n, \   \|A_n \phi\| \leq \|A_{n+1} \phi \| \quad \text{for all} \ \phi \in \D, n \in \N.$$ 
A standard operator algebra (on $\D$) is a *-subalgebra $\mathcal A(\D) \subset \L$ containing the ideal $\mathcal F(\D) \subset \L$ of all finite rank operators on $\D$.\\
Now let $(\mathcal X, \mathcal Y)$ be a dual pair, $\mathcal Y \subset \mathcal X^*$ (the algebraic dual of $\mathcal X$). All linear spaces are supposed to be complex and infinite-dimensional. A linear operator $F \colon \mathcal X \rightarrow \mathcal X$ is said to be of rank one with respect to the dual pair $(\mathcal X, \mathcal Y)$ if $F = x \otimes f$, with nonzero~ $x \in \mathcal X, f \in \mathcal Y$. The linear span 
$$\F(\mathcal X) := \text{span}\{ x \otimes f\colon x \in \mathcal X, f \in \mathcal Y\}$$
is the set of finite-rank operators with respect to the dual pair. If we consider algebras of unbounded operators on $\D$ the ideal $\F(\D) \subset \L$ fits in this concept as follows. Consider the dual pair $(\D, \D^+)$ where $\D^+$ is equal to $\D$ as a set but equipped with the conjugate-linear structure, i.e. the addition in $\D^+$ is the same as in $\D$, but the multiplication by scalars is replaced by the mapping $(\lambda, \phi) \mapsto \bar{\lambda}\phi,~ \lambda \in \C, \phi \in \D$. The pairing is the usual scalar product. Then the rank-one operators from $\F(\D)$ are given by $F = \phi \otimes \psi = \la \psi, \cdot \ra \phi$, i.e. $F \chi = \la \psi, \chi \ra \phi$.\\
Now we are in a position to formulate the results.
\begin{thm}
Let $(\mathcal X, \mathcal Y)$ be a dual pair and let $\F(\mathcal X)$ denote the set of all finite rank operators with respect to the dual pair. Suppose that $\Phi \colon \F(\mathcal X) \rightarrow \F(\mathcal X)$ is an additive surjective mapping preserving operators of rank one in both directions. Then there are either\\
i) a ring automorphism $h \colon \C \rightarrow \C$ and $h$-quasilinear bijective mappings $A \colon \mathcal X \rightarrow \mathcal X,$\\
$ C \colon \mathcal Y \rightarrow \mathcal Y$ such that

$$ \Phi(x \otimes f) = Ax \otimes Cf, \quad x \in \mathcal X, f \in \mathcal Y$$
or \\
ii) a ring automorphism $k \colon \C \rightarrow \C$ and   $k$-quasilinear bijective mappings $A \colon \mathcal X \rightarrow \mathcal Y,$\\
$ C \colon \mathcal Y \rightarrow \mathcal X$ such that

$$\Phi(x \otimes f) = Cf \otimes Ax,  \quad x \in \mathcal X, f \in \mathcal Y$$\end{thm}
\begin{thm}
Let $\D  \subset \H$ be a dense domain in an infinite dimensional complex Hilbert space such that $(\D, t)$ is an (F)-space.   Suppose that $\Phi \colon \F(\D) \rightarrow \F(\D)$ is an additive surjective mapping preserving projections of rank one and their linear spans. Then one of the following four possibilities is valid:\\
i) one has
$$\Phi(T) = ATA^{-1}$$
for all $T \in \F(\D)$, where $A\colon \D \rightarrow \D$ is a bijective either linear or conjugate-linear operator;
ii) one has 
$$\Phi(T) = CT^+C^{-1}$$
for all $T \in \F(\D)$, where $C \colon \D \rightarrow \D$ is a bijective either linear or conjugate-linear operator.\\
If the operators $A,C$ are linear, then they belong to $\L$. 

\end{thm}

\section{Proofs}
The proof of Theorem 1.1 is exactly the proof of Theorem 3.3 in \cite{omladic93}. Let us summarize the main steps of this proof for the readers convenience and for easier reference in course of the proof of Theorem 1.2. We use the following notation. For any nonzero $x \in \mathcal X, f \in \mathcal Y$ denote
$$L_x := \{x \otimes g\colon g \in \mathcal Y\}, \quad R_f := \{u \otimes f\colon u \in \mathcal X\}.$$
\begin{lemma}[cf. \cite{omladic93} Lemmata 2.1 -- 2.4]
Let $\Phi\colon \F(\mathcal X) \rightarrow \F(\mathcal X)$ be an additive and surjective mapping preserving operators of rank one.\\
i) For every $x \in \mathcal X$, either there is a $y \in \mathcal Y$ such that $\Phi(L_x)$ is an additive subgroup of $L_y$ or there is an $f \in \mathcal Y$ such that $\Phi(L_x)$ is an additive subgroup of $R_f$.\\
ii) Let $x$ be a nonzero vector of $\mathcal X$. If $\Phi(L_x)\subset L_y$ then {\rm span }$\Phi(L_x) = L_y$. Also, if $\Phi(L_x) \subset R_f$, then {\rm span }$\Phi(L_x) = R_f$.\\
iii) For every $x,y \in \mathcal X$ we have that either $\Phi(L_x) \subset L_y$ and $\Phi(L_u) \subset L_v$ for appropriate  $u,v \in \mathcal X$, or $\Phi(L_x) \subset R_f$ and $\Phi(L_y) \subset R_g$ for appropriate $f,g \in \mathcal Y$.\\
iv) If for every $x \in \mathcal X$ there exists a $y \in \mathcal X$ such that $\Phi(L_x) \subset L_y$, then for every $f \in \mathcal Y$ there is a $g \in \mathcal Y$ satisfying  $\Phi(R_f) \subset R_g$. Also, if every $\Phi(L_x)$ is a subgroup of an $R_f$ for an $f \in \mathcal Y$, then every $\Phi(R_g)$ is contained in an $L_y$. Consequently, $\Phi(\C x \otimes f) \subset \C\Phi(x \otimes f)$ for every $x \in \mathcal X , f \in \mathcal Y$.
\end{lemma}
\begin{theorem}[cf. \cite{omladic93} Proposition 2.5]
Let $\Phi\colon \F(\mathcal X) \rightarrow \F(\mathcal X)$ be an additive surjective mapping preserving operators of rank one. Then there are either\\
i) a ring homomorphism $h \colon \C \rightarrow \C $ and additive one-to-one mappings $A\colon \mathcal X \rightarrow \mathcal X,$\\
$ C \colon \mathcal Y \rightarrow \mathcal Y$ satisfying $A(\lambda x) = h(\lambda)Ax, C(\lambda f) = h(\lambda)Cf $ for all $\lambda \in \C, x \in \mathcal X, f \in \mathcal Y$ such that
$$\Phi(x \otimes f) = A x \otimes Cf,\quad x \in \mathcal X, f \in \mathcal Y$$

or\\
ii)  a ring homomorphism $k \colon \C \rightarrow \C $ and additive one-to-one mappings $A\colon \mathcal X \rightarrow \mathcal Y,$\\
$ C \colon \mathcal Y \rightarrow \mathcal X$ satisfying $A(\lambda x) = h(\lambda)Ax, C(\lambda f) = h(\lambda)Cf $ for all $\lambda \in \C, x \in \mathcal X, f \in \mathcal Y$ such that
$$\Phi(x \otimes f) = Cf \otimes Ax, \quad x \in \mathcal X, f \in \mathcal Y$$
\end{theorem}
Remark: In Theorem 1.1 the surjectivity of $h, k$ is a consequence of the assumption that $\Phi$ preserves operators of rank one in both directions. Remark that we use the term h-quasilinear only if $h$ is an automorphism and not if it is a homomorphism.\\

To prove Theorem 1.2 it is convenient to formulate Proposition 2.2 for $\F(\D)$.
\begin{theorem}
Let $\D \subset \H$ be a dense domain in a Hilbert space $\H$. Suppose $\Phi \colon \F(\D) \rightarrow \F(\D)$ is an additive surjective mapping preserving operators of rank one. Then there are either\\
i)  a ring homomorphism $h \colon \C \rightarrow \C $ and additive one-to-one mappings $A, C  \colon \D \rightarrow \D$ satisfying $A(\lambda \phi) = h(\lambda)A\phi, C(\lambda \phi) = h(\lambda)C\phi $ for all $\lambda \in \C, \phi \in \D$ such that
$$\Phi(\phi \otimes \psi) = A \phi \otimes C \psi = \la C \psi, \cdot \ra A \phi, \qquad \phi, \psi \in \D$$

or\\
ii) a ring homomorphism $k \colon \C \rightarrow \C $ and additive one-to-one mappings $A, C \colon \D \rightarrow \D$ satisfying $A(\lambda \phi) = k(\lambda)A\phi, C(\lambda \phi) = k(\lambda)C\phi $ for all $\lambda \in \C, \phi \in \D$ such that
$$\Phi(\phi \otimes \psi) = C\psi \otimes A \phi = \la A \phi, \cdot \ra C \psi, \qquad \phi, \psi \in \D$$

\end{theorem} 

The next lemma is also used to prove Theorem 1.2. It is Lemma 4.3 in \cite{omladic93}, and it reflects again a purely algebraic fact. So we omit the proof. We formulate it for more general situations then considered in Theorem 1.2.
\begin{lemma}[\cite{omladic93} Lemma 4.3.] Let $(\mathcal X, \mathcal Y)$ be a dual pair. Assume that $\Phi\colon \F(\mathcal X) \rightarrow \F(\mathcal X)$ is an additive surjective mapping preserving projections of rank one. Then $\Phi$ maps any nilpotent operator of rank one to a nilpotent operator of rank one.
\end{lemma}

{\bf Proof of Theorem 1.2}\\
Lemma 2.4 implies that $\Phi$ preserves operators of rank one. Indeed, if $F = \phi \otimes \psi$ then either $\la \psi, \phi \ra = 0$ (i.e. $F$ is nilpotent and we are done) or $\la \psi, \phi \ra \not= 0$. Then $\lambda F$ and consequently $\Phi(\lambda F) =: P$ is a rank-one projection for an appropriate $\lambda$. Hence $\Phi(F) = \Phi(\frac{1}{\lambda} (\lambda F)) \in \C \Phi(\lambda F)$, i.e. $\Phi(F) = \mu P$. It is easy to see that necessarily $\mu \not= 0$.\\
 So we can apply Proposition 2.3.\\
Let us first consider case i), i.e. $\Phi(\phi \otimes \psi) = A\phi \otimes C \psi$ and let $h, A, C$ be as described in Proposition 2.3 i).\\
 First we show that $A,C$ are bijective mappings which are either both linear or both conjugate-linear. This consists of two steps.\\
Step 1: \\
$$\la C \psi, \chi \ra = h(\la \psi, A^{-1}\chi \ra, \quad \forall \psi \in \D, \chi \in \text{im} A$$
Remember that im $A$ is an additive group and for every $\chi \in \text{im} A$ there is a unique $\phi \in \D$ such that $A\phi = \chi$. If $\la \psi, \phi\ra = 0$ then $\phi \otimes \psi$ is nilpotent of rank one and then so is $\Phi(\phi \otimes \psi) = A\phi \otimes C\psi$, i.e. $\la C\psi, A\phi\ra = 0 = h(\la \psi, \phi \ra) = h(\la \psi, A^{-1}\chi\ra)$.\\
If $\la \psi, \phi \ra \not= 0$ then $[\la \psi, \phi \ra]^{-1} \phi \otimes \psi$ is a rank-one projection and consequently its $\Phi$-image is also a rank-one projection. So
$$\Phi\left([\la \psi, \phi \ra]^{-1} \phi \otimes \psi \right) = h(\la \psi, \phi \ra^{-1}) A\phi \otimes C\psi$$
Using $h(\lambda^{-1}) = h(\lambda)^{-1}$ we get $h(\la \psi, \phi \ra^{-1}) \la C\psi, A\phi \ra = 1$. Therefore
$$\la C \psi, \chi \ra = \la C\psi, A\phi\ra = h(\la \psi, \phi\ra ) = h(\la \psi, A^{-1} \chi\ra).$$
Step 2:\\
$h$ is continuous.  If not, it is unbounded on every neighbourhood of zero (cf. for example \cite{aczel}). Now the proof is almost the same as in  \cite{ti03a} Theorem 3.1. Let $A_n \in \L, n \in \N$ be the operators defining the topology $t$ in $\D$ as described in section 1. Then there are sequences $(\phi_n) \subset \rm{im} A \subset \D, (\psi_n) \subset \D$ such that\\
(i) $ \|A_k \psi_k\| < 2^{-k}$ for all $k \in \N$;\\
(ii) $\| \phi_k\| < 2^{-k}$ for all $k \in \N$;\\
(iii) $\la \psi_i, A^{-1} \phi_k \ra = 0$ for all $i \not= k$;\\
(iv) $| h( \la \psi_n, A^{-1} \phi_n \ra ) | > n + \sum\limits_{i=1}^{n-1} | h(\la \psi_i, A^{-1} \phi_i \ra ) |$.\\ 
These conditions imply that $\psi := \sum\limits_{i=1}^\infty \psi_i \in \D$ and the sequence $(\chi_n)$ with $\chi_n := \sum\limits_{i=1}^n \phi_i$ is $\|\cdot\|$-bounded. Then we have that
$$|\la C\psi, \chi_n\ra| = |h(\la \psi, A^{-1}\chi_n\ra| = |h(\sum\limits_{i=1}^n\la \psi,A^{-1}\phi_i)| = |h(\sum\limits_{i=1}^n\la\psi_i, A^{-1}\phi_i\ra| > n$$
This contradiction proves that $h$ is continuous. Because $h$ is multiplicative it must be either of the form $h(\lambda) = \lambda$ or of the form $h(\lambda) = \bar{\lambda}$ [1, pp. 52 -- 57]. Therefore $A$ and $C$ are either both linear or both conjugate linear. By Lemma 2.1 ii) the mappings $A$ and $C$ must be bijective. Since $\Phi$ preserves projections of rank one and nilpotent operators of rank one, we have for all $\phi, \psi \in \D$
$$\la C \psi, A \phi\ra = \la \psi, \phi \ra\quad A,C - \text{ linear}$$
 or\\
$$\la C \psi, A \phi \ra = \la \phi, \psi \ra \quad A,C - \text{ conjugate-linear}$$
This implies that both operators have adjoints and $\D(C^*) \supset \D, \D(A^*) \supset \D$. Moreover
$$ \la \psi, C^*A\phi \ra = \la \psi, \phi \ra = \la A^*C \psi, \phi\ra$$
which implies $C^*A = I$ on $\D$ or $C = (A^{-1})^+ = (A^+)^{-1}$. If the operators $A,C$ are linear this implies $A,C \in \L$. So we get
$$\Phi(\phi \otimes \psi) = A(\phi \otimes \psi)C^+ = A(\phi \otimes \psi)A^{-1} \quad \text{for all } \phi, \psi \in \D,$$
or
$$\Phi(T) = ATA^{-1} \quad  \text{for all } T \in \F(\D)$$
with $A \colon \D \rightarrow \D$ bijective and either linear or conjugate-linear.\\
This completes the proof in case i).\\
Now let $A,C$ be defined as in case ii) of Proposition 2.3, i.e.
$$\Phi(\phi \otimes \psi) = C \psi \otimes A \phi, \quad \phi, \psi \in \D$$
Quite similar to the previous case it can be shown that $A,C$ are bijective either linear or conjugate-linear operators.Then we get
$$ \Phi(T) = CT^+C^{-1}$$
$C$ either linear or conjugate-linear.  Again, in the linear case we have $C \in \L$.\\
This completes the proof.

\end{document}